\input amstex
\documentstyle{amsppt}
\voffset-,6in
\magnification=\magstep1
\NoBlackBoxes
\topmatter
\title On Lambda-Graph Systems for Subshifts of Subshifts
\endtitle
\author  Wolfgang Krieger
\endauthor
\affil 
Institute for Applied Mathematics\\
University of Heidelberg\\
Im Neuenheimer Feld 294
69120 Heidelberg,
Germany\\
\endaffil
\abstract
Kengo Matsumoto has introduced $\lambda$-graph systems and strong shift
equivalence of $\lambda$-graph systems [Doc.Math.4 (1999), 285-340]. We
associate to a subshift of a subshift a $\lambda$-graph system. The strong
shift equivalence class of the associated  $\lambda$-graph system is an
invariant of subsystem equivalence. Wolfgang Krieger and Kengo Matsumoto have
introduced the $\lambda$-entropy of a $\lambda$-graph system and have shown 
its invariance under strong shift equivalence [Ergod.Th.\&Dynam.Sys. 24 (2004)
1155 - 1172]. A separation entropy of a subshift of a subshift is introduced as the
$\lambda$-entropy of the associated $\lambda$-graph system. 
\endabstract
\endtopmatter



\bigskip
\heading 1. Introduction
\endheading
Let $\Sigma $ be a finite alphabet, and let $S$ denote the shift on $\Sigma
^{\Bbb Z}$,
$$S((x_{i}) _{i \in \Bbb Z}) = ((x_{i+1}) _{i \in \Bbb Z}),\qquad (x_{i}) _{i
\in \Bbb
Z} \in \Bbb Z.
$$
A closed $S$-invariant set $X \subset  \Sigma ^{\Bbb Z} $ with the restriction
of $S$ acting on it, is called a subshift. A finite word is said to
be admissible for a  subshift if it appears in a point of the subshift. A
subshift is uniquely determined by its set of admissible words. A subshift is
said to be of finite type if its admissible words are defined by
excluding finitely many words from appearing as subwords in them. Subshifts are
studied in symbolic dynamics. For an
introduction to symbolic dynamics see \cite {Ki} and \cite {LM}. 

Given subshifts $X \subset
\Sigma ^{\Bbb Z} $ and $\widetilde{X }\subset
\widetilde{\Sigma} ^{\Bbb Z} $, we say that subshifts $Y \subset X$ and
$\widetilde{Y}
\subset \widetilde{X}$ are subsystem equivalent if there exists a topological
conjugacy of $X$ onto $\widetilde {X}$ that carries $Y$ onto $\widetilde {Y}$.
In this
paper we associate to a subshift $Y$ of a subshift $X$ a $\lambda$-graph system
\cite{M1} whose strong shift equivalence class is an invariant of subsystem
equivalence. Consequently the invariants of strong shift equivalence that are
obtained from the $\lambda$-graph systems, e.g. their $\lambda$-entropy \cite
{KM2}, are invariants of subsystem equivalence. 
In section 2 we recall relevant notions and results  and we introduce
terminology and notation.
Section 3 contains a
construction of $\lambda$-graph systems that will be instrumental  in
associating invariantly a $\lambda$-graph system to a subshift of a subshift,
which is done in section 4. In section 5 we consider examples.

\heading 2. Labeled directed graphs and $\lambda$-graph systems
\endheading
We recall from \cite{M1} \cite{KM1}\cite{N} some notions and results concerning labeled
directed graphs, $\lambda$-graph
systems and subshifts. We consider labeled directed graphs $(\Cal V, \Cal E,
\lambda)$ with a set
$\Cal V$ of vertices, a set $\Cal E$ of directed edges and a labeling function
$\lambda$ that assigns to every edge in $\Cal E$ a symbol taken from a finite
alphabet $\Sigma$. (We will suppress the labeling function $\lambda$ in 
the
notation, and also the edge set will not always appear in the notation). A
subset
$\Cal W$ of $\Cal V$ determines a labeled directed
sub-graph of $(\Cal V, \Cal E)  $ with $\Cal W$ as vertex set, in which the
labeled edges in $\Cal E$  are retained that have their initial and final
vertex in $\Cal W$.
A labeled directed graph $(\Cal V , \Cal E )$ is called a Shannon graph if its
labeling is 1-right resolving in the sense that for every vertex $V \in \Cal V$
and every symbol $\sigma \in \Sigma$ there is at most one edge leaving $V$ that
carries the label $\sigma$. We
denote the set of vertices that have an outgoing edge that carries the 
label $\sigma$ by $\Cal V
(\sigma)$, and for $V\in \Cal V (\sigma), \sigma \in \Sigma$ we denote by
$\tau_{\sigma}(V)$ the final vertex of the edge that leaves $V$ and carries
the label $\sigma$. We call $(\tau_{\sigma})_{\sigma \in \Sigma} $ the
transition rules of the Shannon graph. A Shannon graph is determined by its
stransition rules. We say that a Shannon graph such that every vertex has 
an incoming and an outgoing edge presents a subshift $X \subset \Sigma^{\Bbb
Z}$ if the set of
admissible words of the subshift coincides with the set of label sequences of
finite paths in the graph.
By a compact Shannon graph we mean a Shannon graph $(\Cal V, \Cal E)$ where
$\Cal V$ carries a compact topology such that the sets $\Cal 
V(\sigma ), \sigma
\in \Sigma$, are closed (or open) and such that the mappings $V \to
\tau_{\sigma} (V) (V \in \Cal V (\sigma)), \sigma \in \Sigma$, are continuous.
To a compact Shannon graph $ (\Cal V, \Cal E)$ that presents a 
subshift there is associated
the topological Markov chain $M(\Cal V, \Cal E)$,
$$ 
 M(\Cal V, \Cal E) = \bigcap _{i \in \Bbb Z} \{ (V_{i}, x_{i})_{ i\in \Bbb Z}
\in ( \Cal V \times \Sigma)^{\Bbb Z}: V_{i+1} = \tau_{x_{i}}(V_{i}) \}.
$$
The subshift that is presented by a compact Shannon graph $ (\Cal V, \Cal E)$
is identical to the set of label sequences of biinfinite paths on the 
graph,
in other words, $ M(\Cal V, \Cal E) $ projects onto the presented subshift.
For a compact Shannon graph $(\Cal V, \Cal E)$ such that every vertex has an
outgoing edge, let
$$
\tau( \Cal V )= \bigcup_{\sigma \in \Sigma} \tau_{\sigma} (\Cal V(\sigma)),
$$
and set inductively
$$
\tau^{(0)}( \Cal V ) 
= \tau ( \Cal V ), \tau^{(n)}( \Cal V ) = \tau (\tau ^{(n-1)} ( \Cal V 
),\qquad n\in \Bbb N.
$$
Here
$$ 
\tau^{(n)}( \Cal V ) \subset \tau^{(n-1)}( \Cal V ), \qquad n\in \Bbb N,
$$
and the set $\bigcap _{n \in \Bbb N } \tau^{(n)}( \Cal V ) $ determines a
labeled directed sub-Shannon graph of $(\Cal V, \Cal E)$ that presents a
subshift,
and that we denote by  $(\Cal V, \Cal E)^{\circ }$.

A $\lambda$-graph system \cite{M1} is a labeled directed graph in the shape
of a Bratteli diagram with an additional structure. The vertex set of a
$\lambda$-graph system is a disjoint union $\bigcup_{n \in \Bbb Z_{+}} \Cal V_
{n}$ of finite sets, and its edge set is also a disjoint union  $\bigcup_{n
\in \Bbb Z_{+}} \Cal E_ {n}$ of finite sets. Every edge in $\Cal E_ {n}$ has its
initial vertex in $ \Cal V_ {n}$  and its final vertex in $ \Cal V_ 
{n-1}, n \in \Bbb N$. (The edge set  of a $\lambda$-graph system will 
not alway appear 
in the notation.) It is
assumed that all vertices have an incoming edge and that all vertices,
except the vertices in $ \Cal V_ {0}$, have an outgoing edge. The additional
structure that makes the directed labeled Bratteli diagram  into a
$\lambda$-graph system is a shift-like map $\iota :\bigcup_{n \in \Bbb 
N}
\Cal V_ {n}    \to \Cal V   $ such that
$$ 
\iota (   \Cal V_ {n} )  =   \Cal V_ {n-1}   , \qquad     n \in \Bbb N,
$$
that is compatible with the labeling. We say that a $\lambda$-graph system
presents a subshift $X \subset \Sigma ^{\Bbb Z}$ if the set of label
sequences of finite paths in the $\lambda$-graph system coincides with the set
of admissible words of the subshift.
To every Shannon $\lambda$-graph system $\bigcup_{n \in \Bbb Z_{+}} \Cal V_
{n}$ there is associated the compact Shannon
graph
$$
 \{ (V_{n}) _{n \in \Bbb \ N} \in \prod _{n \in \Bbb \ Z_{+}} \Cal V 
_{n}: \iota (V _{n}  )  =  V _{n-1}   , n \in \Bbb N \}
$$
of its $\iota$-orbits, where for $\iota$-orbits $( V _{n})_{n \in 
\Bbb \ Z_{+}} $   and $( W _{n})_{n \in \Bbb \ Z_{+}} $,
$$
\tau_{\sigma}( ( V _{n})_{n \in \Bbb \ Z_{+}} )  = (W _{n})_{n \in \Bbb 
Z_{+} },  
$$
if and only if for $n \in \Bbb N$ 
$\tau_{\sigma}( V  _{n} )   = W _{n-1}  $.
A $\lambda$-graph system is described by its symbolic matrix system $(\Cal M
^{(n, n-1)}, I^{(n, n-1)}  ), n \in  \Bbb N$. Here
$$
\Cal M^{(n, n-1)}       =  (   M^{(n, n-1)}_{c, c^{\prime}})_{c \in \Cal
V_{n},c^{\prime} \in \Cal V_{n-1}      }
$$
gives the symbolic adjacency matrix of the  $\lambda$-graph system:
$$
 M^{(n, n-1)}_{c, c^{\prime}}=\sum _{\{\sigma\in\Sigma:\tau_{\sigma}(c) = 
 c^{\prime}\}} 
 \sigma, \qquad  c \in \Cal
V_{n},c^{\prime} \in \Cal V_{n-1}, 
$$
and
$$
  I^{(n, n-1)}       =  (   I^{(n, n-1)}_{c, c^{\prime}})_{c \in \Cal
V_{n},c^{\prime} \in \Cal V_{n-1}      }
$$
is a 0-1 matrix, that gives the action of the $\iota$-mapping: 
For$ c \in \Cal
V_{n},c^{\prime} \in \Cal V_{n-1}  ,   
n \in \Bbb N$,
$$I^{(n, n-1)}_{c, c^{\prime}} = \cases 1 & \text {if $\iota(c) = 
c^{\prime}$,}\\ 0, &\text{otherwise. } \endcases
$$
The compatibility condition translates into
the commutation relation
$$
\Cal M ^{(n+1, n)} I^{(n, n-1)} = I^{(n+1, n)} \Cal M ^{(n, n-1)}, \qquad n
\in \Bbb N.
$$

Given a finite alphabet $ \Sigma$ we denote by $\iota^{-}$ the operation of
removing the first symbol from a word in $\Sigma^{[1, n]}, n \in \Bbb N$, and by
$\iota^{+}$ the operation of removing the last symbol. $\iota^{-}$ also 
removes the first ysmbol of a sequence in $\Sigma^{\Bbb N}$. We denote by $\Cal
V(\Sigma)$ the set of closed subsets of $\Sigma^{\Bbb N}$ with its 
Hausdorff subset topology, and for
$\sigma \in \Sigma$ we denote
by $\Cal V(\sigma)$ the set of $ V \in \Cal V$ that contain a 
sequence that
begins with $\sigma$. For $\sigma \in\Sigma, V \in \Cal V
_{n}(\Sigma)(\sigma), n \in \Bbb N$, we set 
$$
\tau_{\sigma}(V) = \{ \iota^{-}(v): v \in V, v_{1} = \sigma \}. \tag 1
$$
We have in this way described a
Shannon graph  with vertex set $ 
 \Cal V(\Sigma)$ and transition rules 
$(\tau_{\sigma})_{\sigma \in \Sigma}$.
We denote by $\Cal
V_{n}$ the set of subsets of $\Sigma^{[1, n]}, n \in \Bbb Z_{+}$, and for
$\sigma \in \Sigma$ we denote
by $\Cal V_{n}(\sigma)$ the set of $ V \in \Cal V_{n}$ that contain a word that
begins with $\sigma, n \in \Bbb N$. For $\sigma \in\Sigma, V \in \Cal V
_{n}(\Sigma)(\sigma), n \in \Bbb N$, define
$
\tau_{\sigma}(V) 
$ again by (1).
Using $\iota^{+}$ as the $\iota$-mapping, we have in this way described a
forward separated $\lambda$-graph system with vertex set $ \bigcup_{n \in 
\Bbb Z_{+}} \Cal V_{n}(\Sigma)$ and transition rules 
$(\tau_{\sigma})_{\sigma \in \Sigma}$.
$\Cal
V(\Sigma)$ is the Shannon graph of $\iota$-orbits of the 
$\lambda$-graph system  $ \bigcup_{n \in 
\Bbb Z_{+}} \Cal V_{n}(\Sigma)$.

The forward context $\Gamma ^{+}(V)$ of a vertex $V$ of a Shannon graph is
defined
as the set of sequences in $\Sigma^{\Bbb N}$ that are label sequences of a
semi-infinite path that leaves $V$. A Shannon graph is called forward
separated if distinct vertices have distinct forward contexts. We identify the
vertices  of a forward separated Shannon graph with their forward contexts,
and use then on the vertex set the Hausdorff  subset topology on $\Cal V
(\Sigma)$. A forward separated  Shannon graph whose vertex set  is compact is
a compact Shannon graph. We say that a subset  $\Cal V$ of $ \Cal V (\Sigma)$
is transition complete if for $\sigma \in \Sigma$ and $ V \in \Cal V$ one has
that $\tau_{\sigma} (V) \in \Cal V$. There is a one-to-one correspondence 
between forward separated compact Shannon graphs and the Shannon graphs that
are determined by transition complete closed subsets of $\Cal V (\Sigma)$.
A compact forward separated Shannon graph $\Cal V$ is isomorphic to the
Shannon graph of $\iota $-orbits of a Shannon $\lambda$-graph system that is a
subsystem of $ \bigcup_{n \in \Bbb Z_{+}} \Cal V_{n} (\Sigma)$: one maps a $V
\in \Cal V$ into the $\iota$-orbit $(V_{[1, n]}  )_{ n \in \Bbb Z_{+}  }$.

Let $X \subset \Sigma^{\Bbb Z}$ be a subshift. We denote
$$
	x_{[i,k]} = (x_{j})_{i \leq j \leq k}, \quad  		x \in X,  
			i,k\in \Bbb Z ,  i \leq k,
$$
and
$$
X_{[i,k]} = \{ x_{[i,k]} : x\in X \} . 
$$
We use similar notation  also if indices range in semi-infinite intervals.
Blocks also 
stand for the words
they carry.
We set
$$
 \Gamma ^{+}_{n} (x^{-} ) = \{ a \in X_{[1, n]}:(x^{-}, a) \in X_{( - 
 \infty, n]}\},\qquad n \in \Bbb N,
$$
$$
\Gamma ^{+} (x^{-} ) = \{ x^{+} \in X_{[1, \infty )} : (  x^{-} ,
x^{+}) \in  X  \}   , \qquad  
x^{-} \in  X_{(- \infty, 0 ]},
$$
$$
 \Gamma ^{+}_{n} (a) =  \{ x^{+} \in X_{[1, \infty )} : (  a ,
x^{+}) \in  X_{[-n, \infty )}  \}   , \qquad a \in  X_{[-n.0]}
, n \in \Bbb N.
$$
Every subshift $X \subset \Sigma^{\Bbb Z}$ is presented by its canonical
$\lambda$-graph system with vertex set $\bigcup_{n \in \Bbb Z_{+}} \Cal
V_{n}(X)$,
$$
 \Cal V_{n}(X) = \{ \Gamma^{+}_{n} (x^{-} ) :   x^{-} \in X_{( - \infty, 
 0]}\},
\qquad n \in \Bbb Z_{+},
$$
and by its word  $\lambda$-graph system with vertex set $\bigcup_{n \in \Bbb
Z_{+}} X_{[1, n]}$( or $\bigcup_{n \in \Bbb Z_{+}}  \{\{ x_{[1, n]}\}:
 x \in X  \} $ if one insists that it be a sub-$\lambda$-graph system of    
$\bigcup_{n \in \Bbb Z_{+}} \Cal V_{n}(\Sigma)$). The $\lambda$-graph system
$\bigcup_{n \in \Bbb Z_{+}} \Cal V_{n}(X)$
accompanies the sub-Shannon graph of $\Cal V (\Sigma)$ whose vertex set is the
closure of the set $\{ \Gamma^{+}( x^{-}  ): x^{-} \in X_{( - \infty, 
0]}\} $,
and the $\lambda$-graph system $\bigcup_{n \in \Bbb
Z_{+}} X_{[1, n]}$ accompanies the sub-Shannon graph of $\Cal V
(\Sigma)$ whose vertex set is $X_{[1, \infty)} $ (or $\{\{ x^{+}\}: 
x^{+} \in X_{[1, \infty)} \}$).  

We recall the notion of a bipartite subshift. Let $\Delta $ and $\widetilde
{\Delta}$ be finite disjoint alphabets, and let $Y \subset ( \Delta   \cup
\widetilde {\Delta} )^{\Bbb Z}$ be a subshift. $Y$ is called bipartite if the
admissible words of length two of $Y$ are contained in $\Delta  \widetilde
{\Delta} \cup \widetilde {\Delta} \Delta $. If $Y$ is bipartite then 
$S_{Y}^{2}$
leaves the sets 
$$
X = \{ (y_{i})_{i\in \Bbb Z} \in Y: y_{0} \in  \Delta     \},
$$
and
$$
\widetilde{X }= \{ (y_{i})_{i\in \Bbb Z} \in Y: y_{0} \in  \widetilde {\Delta}
   \},
$$
invariant. Let $S$ resp $\widetilde{S}$ denote the restriction of  
$S_{Y}^{2}$  to $X$
resp. to $ \widetilde {X}$. $(X,S)$ and  $(\widetilde {X},\widetilde {S})$ are
topologically conjugate: a
topological conjugacy of $X$ onto $\widetilde{X}$ is given by the restriction of

$S_{Y}$ to
$X$. Denote the set of words in $  \Delta  \widetilde {\Delta}   $  resp. in $
 \widetilde {\Delta} \Delta $ that are admissible for $Y$ by $\Sigma$ resp. by $
 \widetilde{ \Sigma}     $.
One has $X \subset \Sigma ^{\Bbb Z},    
\widetilde{X} \subset \widetilde { \Sigma} ^{\Bbb Z}$, and one has the
injections
$$
\varphi : \Sigma \hookrightarrow \Delta \widetilde{\Delta}, \quad \widetilde
\varphi: 
\widetilde \Sigma \hookrightarrow \widetilde\Delta \Delta.
$$
By applying $\varphi$ and $\tilde \varphi $ symbol by symbol one extends their
domain of 
definition to finite words and right-infinite sequences.
$\varphi$ and $\widetilde \varphi $ satisfy the relation
$$
\iota^{-}( \iota^{+}    ( \varphi( X_{[1,2]}))) = \widetilde{\varphi}
(\widetilde 
\Sigma), \tag 2
$$
and are called specifications. 
Conversely, let $X \subset \Sigma^{\Bbb Z}$ and $ \widetilde{ X} \subset
\widetilde{\Sigma}^{\Bbb Z}   $ be subshifts, and let $\Delta$ and $\widetilde{
\Delta}     $ be disjoint finite alphabets and let there be given
injections
$$
\varphi : \Sigma \hookrightarrow \Delta \widetilde{\Delta}, \quad 
\widetilde \varphi :\widetilde \Sigma \hookrightarrow
\widetilde\Delta\Delta , 
$$
that are specifications, that is, they satisfy (2). Then
$$
\varphi (X) \cup   \widetilde{\varphi}(\widetilde{X})    \subset \{ \Delta 
\widetilde{\Delta }   \cup 
\widetilde{\Delta } \Delta  \}^{\Bbb Z}
$$
is a bipartite subshift and $X$ and   $ \widetilde{ X} $ are topologically
conjugate. (1) implies that also 
$$
\iota^{-}(\iota^{+}(\varphi (\widetilde X_{[1,2]}))) = \varphi (\Sigma)
$$
and one has a 2-block code given by
$$
a \to \widetilde{\varphi}^{-1}(\iota^{-}(\iota^{+}(\varphi (a)))), \qquad a
\in X_{[1,2]},
$$
that implements a topological conjugacy of $X$ onto $ \tilde{X}$ with the
inverse
given by the 2-block map
$$  
\widetilde{a }\to \varphi^{-1}(\iota^{-}(\iota^{+}(\widetilde{\varphi} 
(\widetilde{a })))), 
\qquad\widetilde{a} \in\widetilde{
X}_{[1,2]}.
$$
Topological conjugacies that arise in this way are called bipartite codings. If
a bipartite coding exists between two subshift presentations then these are
called bipartitely related. 
According to a theorem of Nasu \cite{N} subshifts $X \subset \Sigma 
^{\Bbb Z}     $ and $\widetilde{X} \subset \widetilde{\Sigma} 
^{\Bbb Z}     $are topologically conjugate 
if and only if there is a chain $ X(k) \subset 
\Sigma(k) 
^{\Bbb Z}     , 0 \leq k \leq 
K, K \in \Bbb N$, of subshifts,    $ X(0) = X,  X(K) = \widetilde{X}  $ such
that 
 $ X(k) \subset 
\Sigma(k) 
^{\Bbb Z}     $ and $  X(k+1) \subset 
\Sigma(k+1) 
^{\Bbb Z}     $ are bipartitely related, $ 0 \leq k <
K $.
 A 1-step strong shift equivalence between a symbolic matrix system $(\Cal M
^{(n, n-1)},  I^{(n, n-1)}   )_{n \in \Bbb N}$ with symbol set $\Sigma$ and a
symbolic matrix system $(\widetilde{\Cal M} ^{(n, n-1)},\widetilde{  I}^{(n,
n-1)}   )_{n \in
\Bbb N}    $ with symbol set $\widetilde{\Sigma}$ is given by specifications
$$
\varphi : \Sigma \hookrightarrow \Delta \widetilde{\Delta}, \quad \widetilde
\varphi :\widetilde \Sigma \hookrightarrow
\widetilde\Delta\Delta,
$$
together with symbolic matrices
$$\align
&\Cal K ^{(n, n-1)} = (K_{C, \widetilde{C}})_{ C \in \Cal V_{n}, \widetilde {C
}\in
\widetilde{\Cal V}_{n-1}}    ,\\
&\widetilde{\Cal K}^{(n, n-1)} =(K_{\widetilde{C}, C})_{\tilde{ C 
}\in \widetilde{\Cal V}_{n},  C
\in \Cal V_{n-1} }   , \qquad n \in \Bbb N,
\endalign
$$
such that
$$\align
\Cal K ^{(n+1, n)} \widetilde{\Cal K}^{(n, n-1)} & = \varphi (\Cal M^{(n+1, n)}
I^{(n,n-1)}), \\ 
\widetilde{\Cal K}^{(n+1, n)} \Cal K ^{(n, n-1)} &= \widetilde{\varphi}
(\widetilde{\Cal
M}^{(n+1, n)} I^{(n. n-1)}),\\
\Cal K^{(n+1, n )} \widetilde \varphi ( \widetilde {\Cal M }^{(n, n -1 )})
&=\varphi (\Cal
M ^{(n+1, n)} )\Cal K ^{(n, n-1)},\\
\widetilde{\Cal K }^{(n+1, n ) } \varphi ( \Cal M ^{(n, n -1 )})
&=\tilde{\varphi} (\widetilde{\Cal M} ^{(n+1, n)}) \widetilde{\Cal K} ^{(n, 
n-1)},\\
\Cal K ^{(n+1, n )} \widetilde I  ^{(n, n -1 )} &=I^{(n+1, n)} \Cal K ^{(n, 
n-1)},\\
\widetilde{\Cal K }^{(n+1, n ) } I ^{(n, n -1 )} &=\widetilde{I} ^{(n+1, n)}
\widetilde{\Cal
K} ^{(n, n-1)},\qquad n \in \Bbb N, 
\endalign$$
where the specifications act componentwise on the matrices.
The specifications that come with a 1-step strong shift equivalence between
Shannon $\lambda$-graph systems induce a bipartite coding between the
subshifts that are presented by the Shannon $\lambda$-graph systems.
If subshifts $X \subset \Sigma^{\Bbb Z} $ and $\widetilde{ X }\subset
\widetilde{\Sigma}^{\Bbb Z} $ that are bipartitely related  by specifications
$$
\varphi : \Sigma \hookrightarrow \Delta \widetilde{\Delta}, \quad 
\widetilde {\varphi} :\widetilde \Sigma \hookrightarrow
\widetilde\Delta\Delta ,
$$
then their canonical $\lambda$-graph systems
$\bigcup_{n \in \Bbb Z_{+}} \Cal V_{n}(X)$ and $\bigcup_{n \in \Bbb 
Z_{+}}
\Cal V_{n}(X)$ are strong shift equivalent in one step: A 1-step strong shift
equivalence between their symbolic matrix system is given by the 
symbolic matrices
$$\align
&\Cal K^{(n, n-1)} = (K^{(n, n-1)}_{C,  \widetilde{C}})_{C \in \Cal V_{n} (X),
\widetilde{ C} \in \Cal V_{n-1} (\widetilde{X})}   ,\\
&\widetilde{\Cal K}^{(n, n-1)} = (\widetilde{K}^{(n, n-1)}_{ \widetilde{C},
C})_{\widetilde{C} \in \Cal V_{n} (\widetilde{X}),
\ C \in \Cal V_{n-1} (X)},\qquad  n \in \Bbb N, 
\endalign$$
where 
$$
K^{(n, n-1)}_{C,  \widetilde{C}} = \sum_{\{ \delta \in \iota^{+}(\varphi
(\Sigma ))
: \widetilde{\varphi}^{-1}(\iota^{+}
(\tau_{\delta}(\varphi (C)))) = \widetilde{C}\}} \delta,\qquad C \in \Cal
V_{n} (X), \widetilde{C }
\in \Cal V_{n-1}
(\widetilde{X})     , n \in \Bbb N,
$$
with the symmetric expression for the symbolic matrices $\widetilde{\Cal K
}^{(n, n-1)},
n \in \Bbb N$. The $\lambda$-graph systems $\bigcup _{ n \in \Bbb Z  _{+}}
X_{[1,n]}$ and $\bigcup _{ n \in \Bbb Z  _{+}}\widetilde{ X}_{[1,n]}$ are also
strong shift equivalent in one step: A 1-step strong shift equivalence between
their symbolic matrix systems is given by the symbolic matrices
$$\align
&\Cal K^{(n, n-1)} = ( K^{(n, n-1)} _{c, \tilde{c}}) _{c \in X_{[1, 
n]}, \tilde{c} \in \widetilde{X}_{[1, n-1]}} \\
&\widetilde{\Cal K}^{(n, n-1)} =( \tilde{K}^{(n, n-1)} _{\tilde {c}, c}) _
{\tilde{c} \in \widetilde {X}_{[1, n]}, c \in X_{[1, n-1]}}, \qquad n 
\in \Bbb N,
\endalign$$
where for $c \in X_{[1, 
n]}, \widetilde{c} \in \widetilde{X}_{[1, n-1]}, n \in \Bbb N$,
$$
K^{(n, n-1)}_{c, \tilde{c}} = \cases  \delta  , &\text  
{if $\iota^{+}(\varphi(c_{1})) = \delta, 
\widetilde{\varphi}^{-1}(\iota^{+}(\tau_{\delta}(\varphi(c)))) =
\widetilde{c}$,}     \\
    0 , &\text { otherwise.}\endcases
$$
 with the symmetric expression
for $ \widetilde{\Cal K}^{(n, n-1)},n \in \Bbb N $.  
\heading 3. Pair $\lambda$-graph systems
\endheading

Given Shannon graphs $(\Cal V, \Cal E)$ and $(\Cal W, \Cal F)$, both  with
label set
$\Sigma$, we say that $(\Cal V, \Cal E)$ is subordinate to $(\Cal W, \Cal F)$
if for all $V\in \Cal V$ there exists a $W\in \Cal W$ such that $\Gamma ^{+}
(V) \subset \Gamma ^{+} (W)$. If $(\Cal V, \Cal E)$ is subordinte to $(\Cal
W, \Cal F)$ then we form a pair Shannon graph $[(\Cal V, \Cal E),(\Cal W, \Cal
F) ]$ with vertex set the set of pairs of vertices $(V, W) \in \Cal V \times
\Cal W$ such that $\Gamma ^{+} (V) \subset \Gamma ^{+} (W) $, where the
edges that leave the vertex $(V, W)$ are the pairs that consist of an edge in
$\Cal E$ that leaves $V$ carrying a label $\sigma$ and the edge in $\Cal F$
that leaves $W$ carrying the same label $\sigma$. Given a $\lambda$-graph system
$( \bigcup_{n \in \Bbb Z_{+}}  \Cal V_{n}, \bigcup_{n \in \Bbb Z_{+}}  \Cal
E_{n})$ 
that is subordinate to a
$\lambda$-graph system $(\bigcup _{n \in \Bbb Z_{+}} \Cal W_{n}, \bigcup_{n \in
\Bbb Z_{+}}  \Cal F_{n})$ as a Shannon
graph one turns the pair Shannon graph $[(\bigcup_{n \in \Bbb Z_{+}}  \Cal
V_{n},
\bigcup_{n \in \Bbb Z_{+}}  \Cal
E_{n}),(\bigcup_{n \in \Bbb Z_{+}}  \Cal W_{n}, \bigcup_{n \in \Bbb Z_{+}} 
\Cal F_{n}) ]$ into a pair $\lambda$-graph
system by having the $\iota$-mappings of its constituent $\lambda$-graph systems
acting on the components of its vertices.

Let $X \subset \Sigma^{\Bbb Z}$ be a subshift and let $\Cal V$ be a
forward separated Shannon graph that presents $X$. The Shannon word graph of
the subshift $X$ is subordintate to $\Cal V$ and the word $\lambda$-graph
system $\bigcup_{n \in \Bbb Z_{+}} X_{[1, n ]} $ is subordinate to the 
$\lambda$-graph system $\bigcup_{n \in \Bbb Z_{+}} \Cal V_{[1, n ]} $. One has
here 
$$
([X_{[1, \infty)}, \Cal V])^{\circ}= [X_{[1, \infty)}, \Cal V]
$$
and
$$
([\bigcup_{n \in \Bbb Z_{+}} X_{[1, n]}, \bigcup_{n \in \Bbb Z_{+}} \Cal V_{[1,
n]}])^{\circ}= [\bigcup_{n \in \Bbb Z_{+}} X_{[1, n]}, \bigcup_{n \in \Bbb
Z_{+} }\Cal V_{[1, n]}]
$$
In these pair Shannon graphs resp. pair $\lambda$-graph systems every vertex
has precisely one edge leaving it and this edge carries as label the first
symbol of the word that is the word component of the vertex. We 
denote the vertex set of $[\bigcup_{n \in \Bbb Z_{+}} X_{[1, n]}, \bigcup_{n
\in \Bbb
Z_{+} }\Cal V_{[1, n]}]$ by $ \bigcup_{n \in \Bbb Z_{+}} \widehat{\Cal 
V}_{n} $.
\proclaim {Theorem 3.1}
Let $X \subset \Sigma^{\Bbb Z }$, and $ \tilde{X} \subset \widetilde{X} $
be subshifts that are bipartitely related by specifications
$$
\varphi : \Sigma \hookrightarrow \Delta \widetilde{\Delta},
\quad \widetilde \varphi :\widetilde \Sigma \hookrightarrow
\widetilde\Delta\Delta . 
$$
Let $\Cal V$ be a forward separated Shannon graph that presents $X$, and let
$$
\widetilde{\Cal V} = \{\widetilde{\varphi} ^{-1}(\tau_{\delta}(
(\varphi (\Cal V ))):V \in \Cal V, \delta \in \Delta\} .
$$
Then the $\lambda$-graph systems $
 \bigcup_{n \in \Bbb Z_{+}} \widehat{\Cal V}_{n} $ and $  \bigcup_{n \in \Bbb 
 Z_{+}} \widehat{\widetilde{\Cal V}}_{n}
$ are strong shift equivalent in one step.
\endproclaim
\demo{Proof}Here
$$
\widetilde{\Cal V}_{[1, n]} = \{\widetilde{\varphi} ^{-1}(
\iota^{+}(\tau_{\delta}(
(\varphi (\Cal V )))):
V \in(\Cal V_{[1, n + 1]}, \delta \in \Delta)\}, \qquad n \in \Bbb N.
$$
A 1-step strong shift equivalence between the $\lambda$-graph systems $
 \bigcup_{n \in \Bbb Z_{+}} \widehat{\Cal V}_{n} $ and $  \bigcup_{n \in \Bbb 
 Z_{+}} \widehat{\widetilde{\Cal V}}_{n}$ is
implemented by the symbolic matrices
$$\align
&\Cal K^{(n, n-1)} = 
( K^{(n, n-1)}_{(c,C),( \tilde{c}, \widetilde{C}})
_{(c, C) \in \widehat {\Cal V}_{n} , (\tilde{c},\widetilde{C}) \in
\widehat{\widetilde{ \Cal  V}} _{n-1}},\\
&\widetilde{\Cal K}^{(n, n-1)} =
( K^{(n, n-1)} 
 _{(\tilde{c},\widetilde{C}),( c, C)})
   _{(\tilde{c},\widetilde{C}) \in \widehat{\widetilde{ \Cal V}} _{n},(c, C)
\in
  \widehat {\Cal V}_{n-1}}, \qquad n \in \Bbb N,
\endalign$$
where for $(c, C)
\in \widehat {\Cal V}_{n} , (\widetilde{c},\widetilde{C}) \in
\widehat{\widetilde
{\Cal V}} 
_{n-1}, n \in \Bbb N$,
$$  K^{(n, n-1)} _{(c,C),( \tilde{c}, \widetilde{C})} = \cases \delta  &\text 
{ if $\iota(\varphi(c_{1})) = \delta ,
 \iota^{+}(\tau_{\delta}(\varphi(c))) =\widetilde{c} ,
\widetilde{\varphi}(\iota^{+}\tau_{\delta}(\varphi(C))) = \widetilde {C}$,} \\
0,&\text {otherwise,}
\endcases
$$
with
the symmetric expression
for $ \widetilde{\Cal K}^{(n, n-1)},n \in \Bbb N $. 
\qed
\enddemo
The $\lambda$-entropy  $h_{\lambda}$ of a $\lambda$-graph system was 
introduced in \cite{KM2} as the upper asymptotic growth rate of the 
number of vertices at the n-th level of the $\lambda$-graph system. 
The volume entropy of a  $\lambda$-graph system is 
defined as the upper asymptotic growth rate of the 
number of paths in the $\lambda$-graph system from any
vertex at its $n$-th level to its vertex at level $0$ (see \cite {M2}).
Inspection shows that there is a one-to-one correspondence between the set 
$\{(a, C) \in X_{[1, n ]} \times \Cal V_{[1, n ]}: a \in C\} $ of
vertices that
one finds at the $n$-th level of the $\lambda$-graph system  $[ \bigcup_{n \in
\Bbb Z_{+}} X_{[ 1,n ]}, \bigcup_{n \in \Bbb Z_{+}} \Cal V_{[1, n ]} ]$ and the
set of paths in the $\lambda$-graph system  $  \bigcup_{n \in
\Bbb Z_{+}} \Cal
V_{[1, n ]} $ from any vertex at its $n$-th level to its vertex at level $0$,
$n \in \Bbb N$. The $\lambda$-entropy  of the $\lambda$-graph system
 $[ \bigcup_{n \in \Bbb Z_{+}} X_{[ 1,n ]},  \bigcup_{n \in
\Bbb Z_{+}}\Cal
V_{[1, n ]} ]$ is therefore equal to the volume entropy 
 of the $\lambda$-graph system $  \bigcup_{n \in
\Bbb Z_{+}}
\Cal V_{[1, n ]} $.

\heading 4. Associating a $\lambda$-graph system to a subshift of a subshift
\endheading
Consider a subshift $X \subset \Sigma^{\Bbb Z}$ and a subshift $Y \subset X$.
The Shannon graph $\Cal V (Y)$ and its $\lambda$-graph system $ \bigcup_{n \in
\Bbb Z_{+}}\Cal V _{n} (Y)$ are subordinate to the Shannon graph  $\Cal V
(X)$ and its $\lambda$-graph system $ \bigcup_{n \in
\Bbb Z_{+}}\Cal V _{n}
(X)$. $\Cal V (Y,X)$  will denote the Shannon graph $[\Cal V(Y),\Cal
V(X)]^{\circ} $. $\Cal V (Y,X)$ is the Shannon graph of $\iota$-orbits of the
$\lambda$-graph system
$$
[\bigcup_{n \in \Bbb Z_{+}} \Cal V _{n} (Y), \bigcup_{n \in
\Bbb Z_{+}} \Cal V
_{n} (X)]^{\circ}
$$
whose vertex set we denote by $ \bigcup_{n \in
\Bbb Z_{+}}\Cal V _{n} (Y,X)$.
Note that
$$
\Cal V_{n}(Y,X) = \Cal V (Y,X)_{[1,n]}, \qquad n \in \Bbb N.
$$
We give an alternate description of the vertices of the $\lambda$-graph system
$\bigcup_{n \in \Bbb Z_{+}} \Cal V _{n} (Y,X)    $.

\proclaim {Theorem 4.1}
For all $n \in \Bbb N$ and for
$$
(C_{Y}, C_{X}) \in \Cal V_{n}(Y) \times \Cal V_{n}(X)
$$
one has
$$
(C_{Y}, C_{X}) \in \Cal V_{n}(Y, X)
$$
precisely if
$$
C_{Y}\subset C_{X},
$$
and if for all $K \in \Bbb N$ there is a $k \geq K,$ together with 
$$
b \in Y_{[- k, 0]}, \quad  a\in X_{[- k, 0]},
$$
such that
$$
b _{[- K, 0]}=  a_{[- K, 0]},
$$
and
$$
b^{(l)} \in Y_{[- k-l, 0]}, \quad  a^{(l)}\in X_{[- k-l, 0]}, \qquad l\in \Bbb 
N,
$$
such that
$$
b^{(l)}_ {[- k, 0]} = b, \quad  a^{(l)}_{[- k, 0]} = a,
$$
and
$$
\Gamma^{+}_{Y,n}(b^{(l)}) = \Gamma^{+}_{Y,n}(b),\quad \Gamma^{+}_{X,n}(a^{(l)})
= 
\Gamma^{+}_{X,n}(a),\qquad  l \in \Bbb N.
$$
\endproclaim
\demo{Proof}
The fiber product of the closed subsystem
$$\{ (x_{i}, V_{i}) _{i \in 
\Bbb Z} \in M(\Cal V ) :  (x_{i}) _{i \in  \Bbb Z} \in Y \}$$
of the topological
Markov chain $M(\Cal V)$ and of the topological Markov chain $M(\Cal V)$ with
respect to the projection onto $Y$ yields the topological Markov chain
$$\align
&M(\Cal V (Y, X)) = \\  &\bigcap_{i \in \Bbb Z}\{ ( V^{Y}_{i},  V^{X}_{i} , y
_{i}  )_{i \in \Bbb Z} \in
(\Cal Y \times \Cal Y \times \Sigma)^{\Bbb Z}: V^{Y}_{i}\subset V^{X}_{i},  
V^{Y}_{i+1}   =  \tau_{y_{i}} ( V^{Y}_{i}),   V^{X}_{i+1}   =  \tau_{y_{i}} (
V^{X}_{i}) \}. 
\endalign$$
An examination of the structure of a point in  $M(\Cal V (Y, X))$  confirms the
theorem.
\qed
\enddemo

For subshifts $Y \subset X \subset \Sigma^{\Bbb Z}$ we call the 
$\lambda$-entropy of the $\lambda$-graph system $\bigcup_{ n \in \Bbb
Z_{+}}\Cal V_{n} (Y, X)$ the separation entropy of $Y$ (as a 
subsystem of $X$). One has that 
$$
h_{\lambda} (\bigcup_{ n \in \Bbb Z_{+}} \Cal V_{n}(Y) ) \leq h_{\lambda}
(\bigcup_{ n \in \Bbb Z_{+}} \Cal V_{n}(Y,X) ).
$$

\proclaim {Theorem 4.2}
Let $X \subset \Sigma ^{\Bbb Z} $ and $ \widetilde{X} \subset
\Sigma ^{\Bbb Z}  $ be subshifts and let $ \varphi : X \to \widetilde{ X} $ be a
topological conjuacy of $X$ onto $\widetilde{X}$. Let $Y \subset X$ be a
subshift,
and let $\widetilde{ Y} = \varphi(Y)$. Then the $\lambda$-graph systems
$
\bigcup_{ n \in \Bbb Z_{+}} \Cal V_{n}(Y,X) 
$
and
$
\bigcup_{ n \in \Bbb Z_{+}} \Cal V _{n}(\widetilde{Y},\widetilde{X})
$
are strongly shift equivalent.
 
\endproclaim
\demo{Proof}
By Nasu's theorem it is enough to consider the  case that the subshifts
$X \subset \Sigma ^{\Bbb Z} $   and  $\widetilde{X} \subset \Sigma ^{\Bbb Z}  $
are bipartitely related by a specification
$$
\varphi : \Sigma \hookrightarrow \Delta \widetilde{\Delta}, \quad \widetilde
\varphi 
:\widetilde \Sigma \hookrightarrow
\widetilde\Delta\Delta . 
$$
Denote by $\Sigma_{Y} ( \widetilde{\Sigma}_{\widetilde{Y}})  $ the set of
symbols in $\Sigma
( \widetilde{ \Sigma})$ that are admissible for $Y( \widetilde{ Y})$, and
denote by
$\varphi _{Y}(\widetilde{\varphi} _{\widetilde{Y}}) $ the restriction of
$\varphi_{Y}(
\varphi_{\widetilde{Y}} )$ to $\Sigma_{Y}  (
\widetilde{\Sigma}_{\widetilde{Y}})  $. As is seen
from Theorem (4.1) a 1-step strong shift equivalence between  the
$\lambda$-graph systems
$
\bigcup_{ n \in \Bbb Z_{+}} \Cal V_{n}(Y,X) 
$
and
$
\bigcup_{ n \in \Bbb Z_{+}} \Cal V _{n}(\widetilde{Y},\widetilde{X})
$
is implemented by the symbolic matrices 
$$\align
&\Cal K^{(n, n-1)} = ( K^{(n, n-1)}_{(C_{Y},C_{X}),
( \widetilde{C}_{ \widetilde{Y}}, \widetilde{C}_{
\widetilde{X}})})_{(C_{Y},C_{X}) \in \Cal 
V_{n}(Y,X), ( \widetilde{C}_{ \widetilde{Y}}, 
\tilde{C}_{ \widetilde{X}}) \in \Cal 
V_{n-1}(\widetilde{Y},\widetilde {X})} \\
&\widetilde{\Cal K}^{(n, n-1)} = ( K^{(n, n-1)}_{( \widetilde{C}_{
\widetilde{Y}}, 
\tilde{C}_{, \widetilde{X}}),(C_{Y},C_{X})
})
_{( \widetilde{C}_{ \widetilde{Y}},\widetilde{C}_{ \widetilde{X}}) \in \Cal 
V_{n}(\widetilde{Y},\widetilde {X}) , (C_{Y},C_{X}) \in \Cal 
V_{n-1}(Y,X) }, \qquad n \in \Bbb N, 
\endalign$$
where 
$$\align
K^{(n, n-1)}_{(C_{Y},C_{X}),
( \widetilde{C}_{ \widetilde{Y}}, \widetilde{C}_{ \widetilde{X}})} =
\sum_{\{ \delta 
\in \iota^{+} (\varphi_{Y}(\Sigma_{Y}): 
\iota^{+}(\tau_{\delta}(\varphi_{Y}(C_{Y})))) = \widetilde{C}_{\widetilde{Y}},
\iota^{+}
(\tau_{\delta}(\varphi_{Y}(C_{X})))) = \widetilde{C}_{\widetilde{X}}\}}   & \delta,\\
(C_{Y},C_{X}) \in \Cal 
V_{n}(Y,X), ( \widetilde{C}_{ \widetilde{Y}}, \widetilde{C}_{ \widetilde{X}})
\in \Cal 
V_{n}&(\widetilde{Y},\widetilde {X}),
\endalign$$
with $\widetilde{\Cal K}^{(n, n-1)}
 $  given by the symmetric expression, $n \in \Bbb N$.
\qed
\enddemo

\heading 5. Examples
\endheading

Subsystem equivalence of subshifts of shifts of finite type was studied 
in \cite{BK}. Here we consider some examples of subshifts of 
cartesian powers of a Dyck shift.
We recall the construction of the Dyck shift $D_{2}$. Consider the inverse
monoid with generators $\alpha ^{-}  ,\alpha ^{+} ,\beta ^{-}  ,\beta ^{+}  $
and relations
$$
\alpha ^{-} \alpha ^{+} = \beta ^{-} \beta ^{+} = \bold 1 , 
\alpha ^{-}
 \beta ^{+}  = \beta ^{-}
 \alpha ^{+} = 0.\tag 1
$$
$D_{2}$ is a subshift of $\{\alpha ^{-}  ,\alpha ^{+} ,\beta ^{-}  
,\beta ^{+}   \} ^{ \Bbb Z }$. A point $ (x_{i} )_{{ i \in \Bbb Z}} \in
\{\alpha ^{-}  ,\alpha ^{+} ,\beta ^{-}  
,\beta ^{+}   \} ^{\Bbb Z} $, is in $D_{2}$   precisely if  
$$
\prod_{-I \leq i \leq I} x_{i} \neq 0, \qquad I \in \Bbb N. \tag 2
$$
Let $S_{2}$ denote the shift on $ \{ 0 , 1 \}^{\Bbb Z}$. Setting 
$$
\Phi ^{-} (0)= \alpha ^{-},\Phi ^{-} (1) =
\beta ^{-},  \Phi ^{+} (0)= \alpha ^{+},\Phi ^{+} (1)= \beta ^{+},   
$$
and denoting the identity map on $ \{ \alpha ^{-},
\beta ^{-},  \alpha ^{+}, \beta ^{+} \}$ by $\Phi $, define $ Y^{-}$(
$Y^{+} $) to be the subshift  of $ D_{2} \times  D_{2} $ that is the image of
the
embedding of $ D_{2} \times  S_{2}$ into  $D_{2} \times  D_{2}$ that is given
by the 1-block map $\Phi \times\Phi ^{-}      
 $( $\Phi \times\Phi ^{+}  $), and define  $Y$ to be the subshift of 
 $D_{2} \times  D_{2} \times D_{2} $   that is the image of the 
 embedding of $D_{2} \times  S_{2} \times S_{2} $ into $D_{2} \times  
 D_{2} \times D_{2} $ that is given by the 1-block map $\Phi \times  
 \Phi^{-} \times \Phi^{+} $.   
It is
$$
h_{\lambda} ( \Cal V (D_{2} \times  D_{2} )) = h_{\lambda} 
( \Cal V (Y^{+},  D_{2} \times  S_{2} )) = \log 2 ,   \tag 3
$$
$$
 h_{\lambda}( \Cal V (Y^{-},  D_{2} \times  D_{2}  ) ) = h_{\lambda}
 ( \Cal V (   D_{2} \times  D_{2} ) )=\log 4, \tag 4
$$
and 
$$
 h_{\lambda} ( \Cal V (D_{2} \times  S_{2} \times S_{2}) ) = \log 2,  
$$
$$
 h_{\lambda}( \Cal V ( Y, D_{2} \times  D_{2} \times D_{2}) )=\log 4,
$$
$$
 h_{\lambda} ( \Cal V ( D_{2} \times  D_{2} \times D_{2}) )  = \log 6.
$$
One notes that $D_{2} \times  D_{2} $ has a time reversal that 
carries $Y^{-}$ onto $Y^{+} $. However (3) and (4) imply that there does 
not exist an
automorphism of $D_{2} \times D_{2} $ that carries  $Y^{-} $ onto 
$Y^{+} $.

Let $ K \in \Bbb N$, and consider the semigroup with generators 
$\alpha^{-},  \alpha^{+},\beta^{-},\beta^{+}, $ together
with generators $ \gamma^{-}(k),  \gamma^{+}(k), 1 \leq k \leq K$, with the
relations (1)
and also the relations
$$
 \gamma ^{-}(k) \gamma ^{+}(k) = \bold 1  \qquad 1 \leq k \leq K,   
$$
$$
 \gamma ^{-}(l) \gamma ^{+}(m) = 0, \qquad1 \leq l, m \leq K, l \neq m. 
$$
and the relations
$$
\alpha ^{-}  \gamma ^{+}(k)   =\beta ^{-}  \gamma ^{+}(k) = 0, \ 
\gamma ^{-}(k) \alpha ^{+} 
=\gamma ^{-}(k) \beta ^{+}  =1, \qquad 1 \leq k \leq K.
$$
(comp.  \cite {Kr}).
Define a subshift $X$ of  $(\{ \alpha^{-},  \alpha^{+},
\beta^{-},\beta^{+} \} \cup \{ \gamma^{-}(k),  \gamma^{+}(k), 1 \leq k \leq K\})
^{\Bbb Z}$  that contains the points $ (x_{i} )_{{ i \in \Bbb Z}} \in (\{
\alpha^{-},  \alpha^{+},
\beta^{-},\beta^{+} \} \cup \{ \gamma^{-}(k),  \gamma^{+}(k), 1 \leq k \leq K\})
^{\Bbb Z}$  such that (2) holds.
Here $D_{2} \subset X$ and one has
$$
h_{ \lambda } ( \Cal V ( D_{2}, X))  = 
\log 2 +\log K.
$$
One has 
$$
h_{\lambda} (X)  = \log (2+K).
$$
Note that here for $K > 2$ the separation entropy of the subystem
exceeds the $\lambda$-entropy of the receiving subshift.

\bigskip
{\it e-mail}:

krieger{\@}math.uni-heidelberg.de


\Refs

\refstyle{A}
\widestnumber\key{DGSW}

\ref\key BK
\by  M.Boyle and W.Krieger
\paper Automorphisms and subsystems of the shift
\jour  J. reine
angew. Math.
\vol 437
\yr 1993
\pages 13 - 28
\endref

\ref\key Ki
\by B. P. Kitchens
\book Symbolic dynamics
\publ Springer-Verlag
\publaddr Berlin, Heidelberg and New York
\yr 1998
\endref

\ref\key KM1
\by W.Krieger and K.Matsumoto
\paper Shannon graphs, subshifts and lambda-graph systems
\jour J.Math.Soc.Japan
\vol 54
\yr 2002
\pages 877-899                      
\endref

\ref\key KM2
\by  W.Krieger and K.Matsumoto
\paper A class of invariants of the topological conjugacy of subshifts
\jour  Ergod.
Th.\&Dynam.Sys.
\vol 24
\yr 2004
\pages 1155 - 1172
\endref

\ref\key LM
\by D.Lind and B.Marcus
\book An introduction to symbolic dynamics and coding
\publ Cambridge University Press
\publaddr Cambridge
\yr 1995
\endref

\ref\key M1
\by K.Matsumoto
\paper Presentations of subshifts and there topological conjugacy invariants
\jour Doc.Math,
\vol 4
\yr 1999
\pages  285-340                     
\endref

\ref\key M2
\by K.Matsumoto
\paper Topological entropy in C*-algebras associated with 
lambda-graph systems
\jour  Ergod.
Th.\&Dynam.Sys.
\yr 2005
\pages   1935-1950 
\endref

\ref\key N
\by M.Nasu
\paper Topological conjugacy for sofic shifts
\jour Ergod.Th.\ Dynam.Sys.
\vol 6
\yr 1986
\pages   265-280                   
\endref
\bye